\documentclass[a4paper,12pt]{article}
\usepackage{amssymb,graphics}

\newtheorem{lemma}{Lemma}
\newtheorem{theorem}{Theorem}
\newtheorem{remark}{Remark}
\newtheorem{example}{Example}

\begin{document}

\author{Octavian G. Mustafa\thanks{Correspondence address: Str. Tudor Vladimirescu, Nr. 26, 200534 Craiova, Dolj, Romania}\\
\textit{Faculty of Mathematics, D.A.L.,}\\
\textit{University of Craiova, Romania}\\
\textit{e-mail: octaviangenghiz@yahoo.com}}
\title{On the oscillatory integration of some ordinary differential equations}
\date{}
\maketitle

\noindent\textbf{Abstract} Conditions are given for a class of nonlinear ordinary differential equations $x^{\prime\prime}+a(t)w(x)=0$, $t\geq
t_0\geq1$, which includes the linear equation to possess solutions $x(t)$ with prescribed oblique asymptote that have an oscillatory
pseudo-wronskian $x^{\prime}(t)-\frac{x(t)}{t}$.

\vspace{0.2in} \noindent\textbf{2000 Mathematics Subject Classification:} 34A30 (Primary), 34E05, 34K25


\section{Introduction}

\hspace{0.2in}A certain interest has been shown recently in studying the existence of bounded and positive solutions to a large class of
elliptic partial differential equations which can be displayed as
\begin{eqnarray}
\Delta u+f(x,u)+g(\vert x\vert)x\cdot\nabla u=0,\qquad x\in G_R,\label{gen_class}
\end{eqnarray}
where $G_R=\{x\in\mathbb{R}^n:\vert x\vert>R\}$ for any $R\geq0$ and $n\geq2$. We would like to mention the contributions \cite{AgarwalMustafa},
\cite{Agarwal_et_al}, \cite{Constantin1996} -- \cite{Djebali_et_al}, \cite{Ehrnstrom,EhrnstromMustafa}, \cite{Hesaaraki} and their references in
this respect.

It has been established, see \cite{Constantin1996,Constantin1997}, that it is sufficient for the functions $f$, $g$ to be H\"{o}lder continuous,
respectively continuously differentiable in order to analyze the asymptotic behavior of the solutions to (\ref{gen_class}) by the comparison
method \cite{GilbargTrudinger}. In fact, given $\zeta>0$, let us assume that there exist a continuous function
$A:[R,+\infty)\rightarrow[0,+\infty)$ and a nondecreasing, continuously differentiable function $W:[0,\zeta]\rightarrow[0,+\infty)$ such that
\begin{eqnarray*}
0\leq f(x,u)\leq A(\vert x\vert)W(u)\qquad\mbox{for all }x\in G_R,\thinspace u\in[0,\zeta]
\end{eqnarray*}
and $W(u)>0$ when $u>0$. Then we are interested in the positive solutions $U=U(\vert x\vert)$ of the elliptic partial differential equation
\begin{eqnarray*}
\Delta U+A(\vert x\vert)W(U)=0,\qquad x\in G_R,
\end{eqnarray*}
for the r\^{o}le of super-solutions to (\ref{gen_class}).

M. Ehrnstr\"{o}m \cite{Ehrnstrom} noticed that, by imposing the restriction
\begin{eqnarray*}
x\cdot\nabla U(x)\leq0,\qquad x\in G_R,
\end{eqnarray*}
upon the super-solutions $U$, an improvement of the conclusions from the literature is achieved for the special subclass of equations
(\ref{gen_class}) where $g$ takes only nonnegative values. Further developments of Ehrnstr\"{o}m's idea are given in
\cite{AgarwalMustafa,Agarwal_et_al,Djebali_et_al,EhrnstromMustafa}.

Translated into the language of ordinary differential equations, the research about $U$ reads as follows: given $c_1$, $c_2\geq0$, find (if any)
a positive solution $x(t)$ of the nonlinear differential equation
\begin{eqnarray}
x^{\prime\prime}+a(t)w(x)=0,\qquad t\geq t_0\geq1,\label{the_eq}
\end{eqnarray}
where the coefficient $a:[t_0,+\infty)\rightarrow\mathbb{R}$ and the nonlinearity $w:\mathbb{R}\rightarrow\mathbb{R}$ are continuous and given
by means of $A$, $W$, such that
\begin{eqnarray}
x(t)=c_1t+c_2+o(1)\qquad\mbox{when }t\rightarrow+\infty\label{the_asymptote}
\end{eqnarray}
and
\begin{eqnarray}
{\cal{W}}(x,t)=\frac{1}{t}\left\vert\begin{array}{cc}
x^{\prime}(t)& 1\\
x(t)& t\\
\end{array}\right\vert =x^{\prime}(t)-\frac{x(t)}{t}<0,\qquad t>t_0.\label{pseudo_wronskian}
\end{eqnarray}
The symbol $o(f)$ for a given functional quantity $f$ has here its standard meaning. In particular, by $o(1)$ we refer to a function of $t$ that
decreases to $0$ as $t$ increases to $+\infty$.

The papers \cite{Agarwal_et_al_0,Agarwal_et_al,MustafaRogovchenko0,Mustafa_Glasg,Mustafa2005} present various properties of the functional
quantity ${\cal{W}}$, which shall be called \textit{pseudo-wronskian} in the sequel. Our aim in this note is to complete their conclusions by
giving some sufficient conditions upon $a$ and $w$ which lead to the existence of a solution $x$ to (\ref{the_eq}) that verifies
(\ref{the_asymptote}) while having an \textit{oscillatory pseudo-wronskian} (this means that there exist the unbounded from above sequences
$(t_n^{\pm})_{n\geq1}$ and $(t_n^{0})_{n\geq1}$ such that $t_{2n-1}^0<t_n^{+}<t_{2n}^0<t_n^{-}<t_{2n+1}^0$ and
${\cal{W}}(t_n^{+})>{\cal{W}}(t_n^{0})=0>{\cal{W}}(t_n^{-})$ for all $n\geq1$). We answer thus to a question raised in \cite[p.
371]{Agarwal_et_al}, see also the comment in \cite[pp. 46--47]{Agarwal_et_al_0}.

\section{The sign of ${\cal{W}}$}

\hspace{0.2in} Let us start the discussion with a simple condition to settle the sign issue of the pseudo-wronskian.

\begin{lemma}\label{sign_pseudo_wronskian}
Given $x\in C^{2}([t_0,+\infty),\mathbb{R})$, suppose that $x^{\prime\prime}(t)\leq0$ for all $t\geq t_0$. Then ${\cal{W}}(x,\cdot)$ can change
from being nonnegative-valued to being negative-valued at most once in $[t_0,+\infty)$. In fact, its set of zeros is an interval (possibly
degenerate).
\end{lemma}

\textbf{Proof.} Notice that
\begin{eqnarray*}
\frac{d^2}{dt^2}[x(t)]=\frac{1}{t}\cdot\frac{d}{dt}[t{\cal{W}}(x,t)],\qquad t\geq t_0.
\end{eqnarray*}

The function $t\mapsto t{\cal{W}}(x,t)$ being nonincreasing, it is clear that, if it has zeros, it has either a unique zero or an interval of
zeros. $\square$

The result has an obvious counterpart.

\begin{lemma}\label{sign_pseudo_wronskian1}
Given $x\in C^{2}([t_0,+\infty),\mathbb{R})$, suppose that $x^{\prime\prime}(t)\geq0$ for all $t\geq t_0$. Then, ${\cal{W}}(x,\cdot)$ can change
from being nonpositive-valued to being positive-valued at most once in $[t_0,+\infty)$. Again, its set of zeros is an interval (possibly reduced
to one point).
\end{lemma}

Consider that $x$ is a positive solution of equation (\ref{the_eq}) in the case where $a(t)\geq0$ in $[t_0,+\infty)$ and $w(u)>0$ for all $u>0$.
Then, we have
\begin{eqnarray*}
\frac{d{\cal{W}}}{dt}=-\frac{{\cal{W}}}{t}-a(t)w(x(t)),\qquad t\geq t_0,
\end{eqnarray*}
which leads to
\begin{eqnarray}
{\cal{W}}(x,t)=\frac{1}{t}\left[t_0{\cal{W}}_0-\int_{t_0}^{t}sa(s)w(x(s))ds\right],\qquad
{\cal{W}}_0={\cal{W}}(x,t_0),\label{integral_eq_pseudo_wronskian}
\end{eqnarray}
throughout $[t_0,+\infty)$ by means of Lagrange's variation of constants formula.

The integrand in (\ref{integral_eq_pseudo_wronskian}) being nonnegative-valued, we regain the conclusion of Lemma \ref{sign_pseudo_wronskian}.
In fact, if $T\in[t_0,+\infty)$ is a zero of ${\cal{W}}(x,\cdot)$ then it is a solution of the equation
\begin{eqnarray}
t_0{\cal{W}}_0=\int_{t_0}^{T}sa(s)w(x(s))ds.\label{estimate_T}
\end{eqnarray}
On the other hand, if the pseudo-wronskian of $x$ is positive-valued throughout $[t_0,+\infty)$ then it is necessary to have
\begin{eqnarray}
(t_0{\cal{W}}_0\geq)\qquad\int_{t_0}^{+\infty}sa(s)w(x(s))ds<+\infty.\label{nec_pseudo_wronskian}
\end{eqnarray}

It has become clear at this point that whenever the equation (\ref{the_eq}) has a positive solution $x$ such that ${\cal{W}}_0\leq0$, the
functional coefficient $a$ is nonnegative-valued and has at most isolated zeros and $w(u)>0$ for all $u>0$, the pseudo-wronskian ${\cal{W}}$
satisfies the restriction (\ref{pseudo_wronskian}). Now, returning to the problem stated in the Introduction, we can evaluate the main
difficulty of the investigation: \textit{if the positive solution $x$ has prescribed asymptotic behavior, see formula (\ref{the_asymptote}) or a
similar development, then we cannot decide upfront whether or not ${\cal{W}}_0\leq0$}. The formula (\ref{estimate_T}) shows that there are also
certain difficulties to estimate the zeros of the pseudo-wronskian.

\section{The behavior of ${\cal{W}}$}

\hspace{0.2in} Let us survey in this section some of the recent results regarding the pseudo-wronskian.

It has been established that its presence in the structure of a nonlinear differential equation
\begin{eqnarray}
x^{\prime\prime}+f(t,x,x^{\prime})=0,\qquad t\geq t_0\geq1,\label{the_gen_eq}
\end{eqnarray}
where the nonlinearity $f:[t_0,+\infty)\times\mathbb{R}^2\rightarrow\mathbb{R}$ is continuous, allows for a remarkable flexibility of the
hypotheses when searching for solutions with the asymptotic development (\ref{the_asymptote}) (or similar).

\begin{theorem}
\emph{(\cite[p. 177]{MustafaRogovchenko0})} Assume that there exist the nonnegative-valued, continuous functions $a(t)$ and $g(s)$ such that
$g(s)>0$ for all $s>0$ and $xg(s)\leq g(x^{1-\alpha}s)$, where $x\geq t_0$ and $s\geq0$, for a certain $\alpha\in(0,1)$. Suppose further that
\begin{eqnarray*}
\vert f(t,x,x^{\prime})\vert\leq a(t)g\left(\left\vert
x^{\prime}-\frac{x}{t}\right\vert\right)\quad\mbox{and}\quad\int_{t_0}^{+\infty}\frac{a(s)}{s^{\alpha}}ds<\int_{c+\vert{\cal{W}}_0\vert
t_{0}^{1-\alpha}}^{+\infty}\frac{du}{g(u)}.
\end{eqnarray*}
Then the solution of equation (\ref{the_gen_eq}) given by (\ref{integral_eq_pseudo_wronskian}) exists throughout $[t_0,+\infty)$ and has the
asymptotic behavior
\begin{eqnarray}
x(t)=c\cdot t+o(t),\qquad x^{\prime}(t)=c+o(1)\qquad\mbox{when }t\rightarrow+\infty\label{estim_funkc}
\end{eqnarray}
for some $c=c(x)\in\mathbb{R}$.
\end{theorem}

To compare this result with the standard conditions in asymptotic integration theory regarding the development (\ref{estim_funkc}), see the
papers \cite{Agarwal_et_al_0,Agarwal_et_al,MustafaRogovchenko} and the monograph \cite{KiguradzeChanturia}.

Another result is concerned with the presence of the pseudo-wronskian in the function space $L^{1}((t_0,+\infty),\mathbb{R})$.

\begin{theorem}\label{th_F}
\emph{(\cite[p. 371]{Agarwal_et_al})} Assume that $f$ does not depend explicitly of $x^{\prime}$ and there exists the continuous function
$F:[t_{0},+\infty)\times [0,+\infty)\rightarrow [0,+\infty)$, which is nondecreasing with respect to the second variable, such that
\begin{eqnarray*}
\vert f(t,x)\vert\leq F\left(t,\frac{\vert x\vert}{t}\right)\mbox{ and }\int_{t_{0}}^{+\infty}t\left[1+\ln
\left(\frac{t}{t_{0}}\right)\right]F\left(t,\vert c\vert+\frac{\varepsilon}{t_{0}}\right)dt<\varepsilon
\end{eqnarray*}
for certain numbers $c\neq0$ and $\varepsilon>0$. Then there exists a solution $x(t)$ of equation (\ref{the_gen_eq}) defined in $[t_0,+\infty)$
such that
\begin{eqnarray*}
x(t)=c\cdot t+o(1)\mbox{ when }t\rightarrow+\infty\qquad\mbox{and}\qquad{\cal{W}}(x,\cdot)\in L^1.
\end{eqnarray*}
\end{theorem}

The effect of perturbations upon the pseudo-wronskian is investigated in the papers \cite{Agarwal_et_al_0,MustafaRogovchenko0,Mustafa_Glasg}.

\begin{theorem}
\emph{(\cite[p. 183]{MustafaRogovchenko0})} Consider the nonlinear differential equation
\begin{eqnarray}
x^{\prime\prime}+f(t,x,x^{\prime})=p(t),\qquad t\geq t_0\geq1,\label{the_gen_eq_perturbed}
\end{eqnarray}
where the functions $f:[t_0,+\infty)\times\mathbb{R}^2\rightarrow\mathbb{R}$ and $p:[t_0,+\infty)\rightarrow\mathbb{R}$ are continuous and
verify the hypotheses
\begin{eqnarray*}
\vert f(t,x,x^{\prime})\vert\leq a(t)\left\vert x^{\prime}-\frac{x}{t}\right\vert,\qquad \int_{t_0}^{+\infty}ta(t)dt<+\infty
\end{eqnarray*}
and
\begin{eqnarray*}
\lim\limits_{t\rightarrow+\infty}\frac{1}{t}\int_{t_0}^{t}sp(s)ds=C\in\mathbb{R}-\{0\}.
\end{eqnarray*}
Then, given $x_0\in\mathbb{R}$, there exists a solution $x(t)$ of equation (\ref{the_gen_eq_perturbed}) defined in $[t_0,+\infty)$ such that
\begin{eqnarray*}
x(t_0)=x_0\qquad\mbox{and}\qquad\lim\limits_{t\rightarrow+\infty}{\cal{W}}(x,t)=C.
\end{eqnarray*}
In particular,
\begin{eqnarray*}
\lim\limits_{t\rightarrow+\infty}\frac{x(t)}{t\ln t}=C.
\end{eqnarray*}
\end{theorem}

A slight modification of the discussion in \cite[Remark 3]{Mustafa_Glasg}, see \cite[p. 47]{Agarwal_et_al_0}, leads to the next result.
\begin{theorem}
Assume that $f$ in (\ref{the_gen_eq_perturbed}) does not depend explicitly of $x^{\prime}$ and there exists the continuous function
$F:[t_{0},+\infty)\times [0,+\infty)\rightarrow [0,+\infty)$, which is nondecreasing with respect to the second variable, such that
\begin{eqnarray*}
\vert f(t,x)\vert\leq F(t,\vert x\vert)\mbox{ and }\int_{t}^{+\infty}sF\left(s,\vert P(s)\vert+\sup\limits_{\tau\geq s}\{q(\tau)\}\right)ds\leq
q(t),\thinspace t\geq t_0,
\end{eqnarray*}
for a certain positive-valued, continuous function $q(t)$ possibly decaying to $0$ as $t\rightarrow+\infty$. Here, $P$ is the twice continuously
differentiable antiderivative of $p$, that is $P^{\prime\prime}(t)=p(t)$ for all $t\geq t_0$. Suppose further that
\begin{eqnarray*}
\limsup\limits_{t\rightarrow+\infty}\left[t\frac{{\cal{W}}(P,t)}{q(t)}\right]>1\qquad\mbox{and}\qquad\liminf\limits_{t\rightarrow+\infty}\left[t\frac{{\cal{W}}(P,t)}{q(t)}\right]<-1.
\end{eqnarray*}
Then equation (\ref{the_gen_eq_perturbed}) has a solution $x(t)$ throughout $[t_0,+\infty)$ such that
\begin{eqnarray*}
x(t)=P(t)+o(1)\qquad\mbox{when }t\rightarrow+\infty
\end{eqnarray*}
and ${\cal{W}}(x,\cdot)$ oscillates.
\end{theorem}

Finally, the presence of the pseudo-wronskian in the structure of a nonlinear differential equation can lead to multiplicity when searching for
solutions with the asymptotic development (\ref{the_asymptote}).

\begin{theorem}\emph{(\cite[Theorem 1]{Mustafa2005})}
Given the numbers $x_0$, $x_1$, $c\in\mathbb{R}$, with $c\neq0$, and $t_0\geq1$ such that $t_0x_1-x_0=c$, consider the Cauchy problem
\begin{eqnarray}
\left\{\begin{array}{ll}x^{\prime\prime}=\frac{1}{t}g(tx^{\prime}-x),\qquad t\geq t_0\geq1,\\
x(t_0)=x_0,\qquad x^{\prime}(t_0)=x_1,\end{array}\right.\label{Cauchyproblem}
\end{eqnarray}
where the function $g:\mathbb{R}\times\mathbb{R}\rightarrow\mathbb{R}$ is continuous, $g(c)=g(3c)=0$ and $g(u)>0$ for all $u\neq c$. Assume
further that
\begin{eqnarray*}
\int_{c+}^{2c}\frac{du}{g(u)}<+\infty\qquad\mbox{and}\qquad\int_{2c}^{(3c)-}\frac{du}{g(u)}=+\infty.
\end{eqnarray*}
Then problem (\ref{Cauchyproblem}) has an infinity of solutions $x(t)$ defined in $[t_0,+\infty)$ and developable as
\begin{eqnarray*}
x(t)=c_{1}t+c_2+o(1)\qquad\mbox{when }t\rightarrow+\infty
\end{eqnarray*}
for some $c_1=c_1(x)$ and $c_2=c_2(x)\in\mathbb{R}$.
\end{theorem}

The asymptotic analysis of certain functional quantities attached to the solutions of equations (\ref{the_eq}), (\ref{the_gen_eq}) and
(\ref{the_gen_eq_perturbed}), as in our case the pseudo-wronskian, might lead to some surprising consequences. Among the functional quantities
that gave the impetus to spectacular developments in the qualitative theory of linear/nonlinear ordinary differential equations we would like to
refer to
\begin{eqnarray*}
{\cal{K}}(x)(t)=x(t)x^{\prime}(t),\qquad t\geq t_0,
\end{eqnarray*}
employed in the theory of \textit{Kneser-solutions}, see the papers \cite{CMV1991,CMV1999} for the linear and respectively the nonlinear case
and the monograph \cite{KiguradzeChanturia}, and
\begin{eqnarray*}
{\cal{HW}}(x)=\int_{t_0}^{+\infty}x(s)w(x(s))ds.
\end{eqnarray*}
The latter quantity is the core of the nonlinear version of \textit{Hermann Weyl's limit-point/limit-circle classification} designed for
equation (\ref{the_eq}), see the well-documented monograph \cite{BDG} and the paper \cite{MustafaRogovchenkoJMAA}.

\section{The negative values of ${\cal{W}}$}

\hspace{0.2in} We shall assume in the sequel that the nonlinearity $w$ of equation (\ref{the_eq}) verifies some of the hypotheses listed below:
\begin{eqnarray}
\vert w(x)-w(y)\vert\leq k\vert x-y\vert,\qquad\mbox{where }k>0,\label{hyp1_w}
\end{eqnarray}
and
\begin{eqnarray}
w(0)=0,\thinspace w(x)>0\mbox{ when }x>0,\qquad \vert w(xy)\vert\leq w(\vert x\vert)w(\vert y\vert)\label{hyp2_w}
\end{eqnarray}
for all $x$, $y\in\mathbb{R}$. We notice that restriction (\ref{hyp2_w}) implies the existence of a majorizing function $F$, as in Theorem
\ref{th_F}, given by the estimates
\begin{eqnarray*}
\vert f(t,x)\vert=\vert a(t)w(x)\vert\leq \vert a(t)\vert\cdot w(t)w\left(\frac{\vert x\vert}{t}\right)=F\left(t,\frac{\vert x\vert}{t}\right).
\end{eqnarray*}

We can now use the paper \cite{MustafaRogovchenko} to recall the main conclusions of an asymptotic integration of equation (\ref{the_eq}). It
has been established that whenever $\int^{+\infty}_{t_0}tw(t)\vert a(t)\vert dt<+\infty$, all the solutions of (\ref{the_eq}) have asymptotes
(\ref{the_asymptote}) and their first derivatives are developable as
\begin{eqnarray}
x^{\prime}(t)=c_1+o\left(t^{-1}\right)\qquad\mbox{when }t\rightarrow+\infty.\label{the_derivative}
\end{eqnarray}
Consequently, ${\cal{W}}(x,t)=-c_2t^{-1}+o(t^{-1})$ for all large $t$'s. In this case (the functional coefficient $a$ has varying sign), when
dealing with the sign of the pseudo-wronskian, of interest would be the subcase where $c_2=0$. Here, the asymptotic development does not even
ensure that ${\cal{W}}$ is eventually negative. Enlarging the family of coefficients to the ones subjected to the restriction
$\int_{t_0}^{+\infty}t^{\varepsilon}w(t)\vert a(t)\vert dt<+\infty$, where $\varepsilon\in[0,1)$, the developments (\ref{the_asymptote}),
(\ref{the_derivative}) become
\begin{eqnarray}
x(t)=ct+o\left(t^{1-\varepsilon}\right),\qquad x^{\prime}(t)=c+o\left(t^{-\varepsilon}\right),\qquad
c\in\mathbb{R},\label{less_precise_development}
\end{eqnarray}
yielding the less precise estimate ${\cal{W}}(x,t)=o(t^{-\varepsilon})$ when $t\rightarrow+\infty$. We have again a lack of precision in the
asymptotic development of ${\cal{W}}(x,\cdot)$ with respect to the sign issue. We also deduce on the basis of (\ref{the_asymptote}),
(\ref{less_precise_development}) that some of the coefficients $a$ in these classes verify (\ref{nec_pseudo_wronskian}), a fact that complicates
the discussion.

The next result establishes the existence of a positive solution to (\ref{the_eq}) subjected to (\ref{pseudo_wronskian}),
(\ref{less_precise_development}) for the largest class of functional coefficients: $\varepsilon=0$. By taking into account Lemmas
\ref{sign_pseudo_wronskian}, \ref{sign_pseudo_wronskian1} and the non-oscillatory character of equation (\ref{the_eq}) when the nonlinearity $w$
verifies (\ref{hyp2_w}), we conclude that \textit{for an investigation within this class of coefficients $a$ of the solutions with oscillatory
pseudo-wronskian it is necessary that $a$ itself oscillates}. Also, when $a$ is non-negative valued we recall that the condition
\begin{eqnarray*}
\int_{t_0}^{+\infty}a(t)dt<+\infty
\end{eqnarray*}
is necessary for the linear case of equation (\ref{the_eq}) to be non-oscillatory, see \cite{Hartman}, while in the case given by
$w(x)=x^{\lambda}$, $x\in\mathbb{R}$, with $\lambda>1$ (such an equation is usually called an \textit{Emden-Fowler equation}, see the monograph
\cite{KiguradzeChanturia}) the condition
\begin{eqnarray}
\int_{t_0}^{+\infty}ta(t)dt=+\infty\label{Atkinson_condition}
\end{eqnarray}
is necessary and sufficient for oscillation, see \cite{Atkinson}. In the case of Emden-Fowler equations with $\lambda\in(0,1)$ and a
continuously differentiable coefficient $a$ such that $a(t)\geq0$ and $a^{\prime}(t)\leq0$ throughout $[t_0,+\infty)$, another result
establishes that equation (\ref{the_eq}) has no oscillatory solutions provided that condition (\ref{Atkinson_condition}) fails, see
\cite{Heidel}.

Regardless of the oscillation of $a$, it is known \cite[p. 360]{Agarwal_et_al} that the linear case of equation (\ref{the_eq}) has bounded and
positive solutions with eventually negative pseudo-wronskian.

\begin{theorem}
Assume that the nonlinearity $w$ verifies hypothesis (\ref{hyp2_w}) and is nondecreasing. Given $c$, $d>0$, suppose that the functional
coefficient $a$ is nonnegative-valued, with eventual isolated zeros, and
\begin{eqnarray*}
\int_{t_0}^{+\infty}w(t)a(t)dt\leq\frac{d}{w(c+d)}.
\end{eqnarray*}
Then, the equation (\ref{the_eq}) has a solution $x$ such that ${\cal{W}}_0=0$,
\begin{eqnarray}
c-d\leq x^{\prime}(t)<\frac{x(t)}{t}\leq c+d\qquad\mbox{for all }t>t_0\label{est_1}
\end{eqnarray}
and
\begin{eqnarray}
\lim\limits_{t\rightarrow+\infty}x^{\prime}(t)=\lim\limits_{t\rightarrow+\infty}\frac{x(t)}{t}=c.\label{est_2}
\end{eqnarray}
\end{theorem}

\textbf{Proof.} We introduce the set $D$ given by
\begin{eqnarray*}
D=\{u\in C([t_0,+\infty),\mathbb{R}):ct\leq u(t)\leq (c+d)t\mbox{ for every }t\geq t_0\}.
\end{eqnarray*}

A partial order on $D$ is provided by the usual pointwise order $"\leq"$, that is, we say that $v_1\leq v_2$ if and only if $v_1(t)\leq v_2(t)$
for all $t\geq t_1$, where $v_1$, $v_2\in D$. It is not hard to see that $(D,\leq)$ is a complete lattice.

For the operator $V:D\rightarrow C([t_0,+\infty),\mathbb{R})$ with the formula
\begin{eqnarray*}
V(u)(t)=t\left\{c+\int_{t}^{+\infty}\frac{1}{s^{2}}\int_{t_0}^{s}\tau a(\tau)w(u(\tau))d\tau ds\right\},\qquad u\in D,\thinspace t\geq t_0,
\end{eqnarray*}
the next estimates hold
\begin{eqnarray*}
c&\leq&\frac{V(u)(t)}{t}= c+\int_{t}^{+\infty}\frac{1}{s^{2}}\int_{t_0}^{s}\tau a(\tau)\cdot w(\tau)w\left(\frac{u(\tau)}{\tau}\right)d\tau ds\\
&\leq& c+\sup\limits_{\xi\in[0,c+d]}\{w(\xi)\}\cdot\int_{t}^{+\infty}\frac{1}{s^{2}}\int_{t_0}^{s}\tau w(\tau)
a(\tau)d\tau ds\\
&=&c+w(c+d)\left[\frac{1}{t}\int_{t_0}^{t}\tau w(\tau)
a(\tau)d\tau+\int_{t}^{+\infty}w(\tau) a(\tau)d\tau\right]\\
&\leq&c+w(c+d)\int_{t_0}^{+\infty}w(\tau) a(\tau)d\tau\leq c+d
\end{eqnarray*}
by means of (\ref{hyp2_w}). These imply that $V(D)\subseteq D$.

Since $c\cdot t\leq V(c\cdot t)$ for all $t\geq t_0$, by applying the Knaster-Tarski fixed point theorem \cite[p. 14]{DG}, we deduce that the
operator $V$ has a fixed point $u_0$ in $D$. This is the pointwise limit of the sequence of functions $(V^n(c\cdot\mbox{Id}_{I}))_{n\geq1}$,
where $V^1=V$, $V^{n+1}=V^n\circ V$ and $I=[t_0,+\infty)$.

We deduce that
\begin{eqnarray*}
u_{0}^{\prime}(t)&=&[V(u_0)]^{\prime}(t)=\frac{u_0(t)}{t}-\frac{1}{t}\int_{t_0}^{t}\tau a(\tau)w(u_0(\tau))d\tau<\frac{u_0(t)}{t},
\end{eqnarray*}
when $t>t_0$, and thus (\ref{est_1}), (\ref{est_2}) hold true.

The proof is complete. $\square$

\section{The oscillatory integration of equation (\ref{the_eq})}

\hspace{0.2in} Let the continuous functional coefficient $a$ with varying sign satisfy the restriction
\begin{eqnarray*}
\int_{t_0}^{+\infty}t^2\vert a(t)\vert dt<+\infty.
\end{eqnarray*}

We call the problem studied in the sequel \textit{an oscillatory (asymptotic) integration} of equation (\ref{the_eq}).

\begin{theorem}\label{main_result}
Assume that $w$ verifies (\ref{hyp1_w}), $w(0)=0$ and there exists $c>0$ such that
\begin{eqnarray}
L_{+}^{c}>0>L_{-}^{c},\label{L_c}
\end{eqnarray}
where
\begin{eqnarray*}
L_{+}^{c}=\limsup\limits_{t\rightarrow+\infty}\frac{t\int_{t}^{+\infty}sw(cs)a(s)ds}{\int_{t}^{+\infty}s^2\vert a(s)\vert ds},\quad
L_{-}^{c}=\liminf\limits_{t\rightarrow+\infty}\frac{t\int_{t}^{+\infty}sw(cs)a(s)ds}{\int_{t}^{+\infty}s^2\vert a(s)\vert ds}.
\end{eqnarray*}
Then the equation (\ref{the_eq}) has a solution $x(t)$ with oscillatory pseudo-wronskian such that
\begin{eqnarray}
x(t)=c\cdot t+o(1)\qquad\mbox{when }t\rightarrow+\infty.\label{the_new_asymptote}
\end{eqnarray}
\end{theorem}

\textbf{Proof.} There exist $\eta>0$ such that $L_{+}^{c}>\eta$, $L_{-}^{c}<-\eta$ and two increasing, unbounded from above sequences
$(t_n)_{n\geq1}$, $(t^n)_{n\geq1}$ of numbers from $(t_0,+\infty)$ such that $t^n\in(t_n,t_{n+1})$ and
\begin{eqnarray}
t_n\int_{t_n}^{+\infty}sw(cs)a(s)ds+k\eta\int_{t_n}^{+\infty}s^2\vert a(s)\vert ds<0\label{t_n1}
\end{eqnarray}
and
\begin{eqnarray}
t^n\int_{t^n}^{+\infty}sw(cs)a(s)ds-k\eta\int_{t^n}^{+\infty}s^2\vert a(s)\vert ds>0\label{t_n2}
\end{eqnarray}
for all $n\geq1$.

Assume further that
\begin{eqnarray*}
\int_{t_0}^{+\infty}\tau^2\vert a(\tau)\vert d\tau\leq\frac{\eta}{k(c+\eta)}
\end{eqnarray*}
and introduce the complete metric space $S=(D,\delta)$ given by
\begin{eqnarray*}
D=\{y\in C([t_0,+\infty),\mathbb{R}):t\vert y(t)\vert\leq\eta\mbox{ for every }t\geq t_0\}
\end{eqnarray*}
and
\begin{eqnarray*}
\delta(y_1,y_2)=\sup\limits_{t\geq t_0}\{t\vert y_1(t)-y_2(t)\vert\},\qquad y_1,\thinspace y_2\in D.
\end{eqnarray*}

For the operator $V:D\rightarrow C([t_0,+\infty),\mathbb{R})$ with the formula
\begin{eqnarray*}
V(y)(t)=\frac{1}{t}\int_{t}^{+\infty}sa(s)w\left(s\left[c-\int_{s}^{+\infty}\frac{y(\tau)}{\tau}d\tau\right]\right)ds,\quad y\in D,\thinspace
t\geq t_0,
\end{eqnarray*}
the next estimates hold (notice that $\vert w(x)\vert\leq k\vert x\vert$ for all $x\in\mathbb{R}$)
\begin{eqnarray}
t\vert V(y)(t)\vert\leq k\int_{t}^{+\infty}s^2\vert a(s)\vert\left[c+\eta\int_{s}^{+\infty}\frac{d\tau}{\tau^2}\right]ds\leq\eta\label{V_eta}
\end{eqnarray}
and
\begin{eqnarray*}
t\vert V(y_2)(t)-V(y_1)(t)\vert&\leq&k\int_{t}^{+\infty}s^2\vert
a(s)\vert\left(\int_{s}^{+\infty}\frac{d\tau}{\tau^2}\right)ds\cdot\delta(y_1,y_2)\\
&\leq&\frac{k}{t_0}\int_{t}^{+\infty}s^2\vert a(s)\vert ds\leq\frac{\eta}{c+\eta}\cdot\delta(y_1,y_2).
\end{eqnarray*}
These imply that $V(D)\subseteq D$ and thus $V:S\rightarrow S$ is a contraction.

From the formula of operator $V$ we notice also that
\begin{eqnarray}
\lim\limits_{t\rightarrow+\infty}tV(y)(t)=0\qquad\mbox{for all }y\in D.\label{est_3}
\end{eqnarray}

Given $y_0\in D$ the unique fixed point of $V$, one of the solutions to (\ref{the_eq}) has the formula
$x_0(t)=t\left[c-\int_{t}^{+\infty}\frac{y_0(s)}{s}ds\right]$ for all $t\geq t_0$. Via (\ref{est_3}) and L'Hospital's rule, we provide also an
asymptotic development for this solution, namely
\begin{eqnarray*}
\lim\limits_{t\rightarrow+\infty}[x_0(t)-c\cdot t]&=&-\lim\limits_{t\rightarrow+\infty}t\int_{t}^{+\infty}\frac{y_0(s)}{s}ds
=-\lim\limits_{t\rightarrow+\infty}ty_0(t)\\
&=&-\lim\limits_{t\rightarrow+\infty}tV(y_0)(t)=0.
\end{eqnarray*}

The estimate
\begin{eqnarray*}
\left\vert ty_0(t)-\int_{t}^{+\infty}sw(cs)a(s)ds\right\vert&\leq&k\int_{t}^{+\infty}s^2\vert a(s)\vert\left[\int_{s}^{+\infty}\frac{\vert
y_0(\tau)\vert}{\tau}d\tau\right]ds\\
&\leq&k\eta\cdot\frac{1}{t}\int_{t}^{+\infty}s^2\vert a(s)\vert ds,\qquad t\geq t_0,
\end{eqnarray*}
accompanied by (\ref{t_n1}), (\ref{t_n2}), leads to
\begin{eqnarray}
y_0(t_n)={\cal{W}}(x_0,t_n)<0\qquad\mbox{and}\qquad y_0(t^n)={\cal{W}}(x_0,t^n)>0.\label{estim_oscil}
\end{eqnarray}

The proof is complete. $\square$

\begin{remark}
\emph{When equation (\ref{the_eq}) is linear, that is $w(x)=x$ for all $x\in\mathbb{R}$, the formula (\ref{L_c}) can be recast as
\begin{eqnarray*}
L_{+}=\limsup\limits_{t\rightarrow+\infty}\frac{t\int_{t}^{+\infty}s^2a(s)ds}{\int_{t}^{+\infty}s^2\vert a(s)\vert
ds}>0>\liminf\limits_{t\rightarrow+\infty}\frac{t\int_{t}^{+\infty}s^2a(s)ds}{\int_{t}^{+\infty}s^2\vert a(s)\vert ds}=L_{-}.
\end{eqnarray*}
We claim that} for all $c\neq0$ there exists a solution $x(t)$ with oscillatory pseudo-wronskian which verifies (\ref{the_new_asymptote}).
\emph{In fact, replace $c$ with $c_0$ in the formulas (\ref{t_n1}), (\ref{t_n2}) for a certain $c_0$ subjected to the inequality
$\min\{L_{+},-L_{-}\}>\frac{\eta}{c_0}$. It is obvious that, when $L_{+}=-L_{-}=+\infty$, formulas (\ref{t_n1}), (\ref{t_n2}) hold for all
$c_0$, $\eta>0$. Given $c\in\mathbb{R}-\{0\}$, there exists $\lambda\neq0$ such that $c=\lambda c_0$. The solution of equation (\ref{the_eq})
that we are looking for has the formula $x=\lambda\cdot x_0$, where $x_{0}(t)=t\left[c_0-\int_{t}^{+\infty}\frac{y_0(s)}{s}ds\right]$ for all
$t\geq t_0$ and $y_0$ is the fixed point of operator $V$ in $D$. Its pseudo-wronskian oscillates as a consequence of the obvious identity
\begin{eqnarray*}
\lambda\cdot{\cal{W}}(x_0,t)={\cal{W}}(x,t),\qquad t\geq t_0.
\end{eqnarray*}
}
\end{remark}

\begin{example}
\emph{An immediate example of functional coefficient $a$ for the problem of linear oscillatory integration is given by $a(t)=t^{-2}e^{-t}\cos
t$, where $t\geq1$.}

\emph{We have
\begin{eqnarray*}
\int_{t}^{+\infty}s^2a(s)ds=\frac{1}{\sqrt{2}}\cos\left(t+\frac{\pi}{4}\right)e^{-t}\quad\mbox{and}\quad\int_{t}^{+\infty}s^2\vert a(s)\vert
ds\leq e^{-t}
\end{eqnarray*}
throughout $[1,+\infty)$ which yields $L_{+}=+\infty$, $L_{-}=-\infty$.}
\end{example}

Sufficient conditions are provided now for an oscillatory pseudo-wronskian to be in $L^{p}((t_0,+\infty),\mathbb{R})$, where $p>0$. Since
$\lim\limits_{t\rightarrow+\infty}{\cal{W}}(x,t)=0$ for any solution $x(t)$ of equation (\ref{the_eq}) with the asymptotic development
(\ref{the_new_asymptote}), (\ref{the_derivative}), we are interested in the case $p\in(0,1)$.

\begin{theorem}\label{oscil_L_p}
Assume that, in the hypotheses of Theorem \ref{main_result}, the coefficient $a$ verifies the condition
\begin{eqnarray}
\int_{t_0}^{+\infty}\left[\frac{t}{\int_{t}^{+\infty}s^2\vert a(s)\vert ds}\right]^{1-p}t^2\vert a(t)\vert dt<+\infty\quad\mbox{for some
}p\in(0,1).\label{hyp_a_L_p}
\end{eqnarray}
Then the equation (\ref{the_eq}) has a solution $x(t)$ with an oscillatory pseudo-wronskian in $L^p$ and the asymptotic expansion
(\ref{the_new_asymptote}).
\end{theorem}

\textbf{Proof.} Recall that $y_0$ is the fixed point of operator $V$. Then, formula (\ref{V_eta}) implies that
\begin{eqnarray*}
\vert y_0(t)\vert\leq k(c+\eta)\cdot\frac{1}{t}\int_{t}^{+\infty}s^2\vert a(s)\vert ds,\qquad t\geq t_0.
\end{eqnarray*}

Via an integration by parts, we have
\begin{eqnarray*}
&&\frac{1}{[k(c+\eta)]^p}\int_{t}^{T}\vert y_{0}(s)\vert^p ds\\
&&\leq\frac{T^{1-p}}{1-p}\left[\int_{T}^{+\infty}s^2\vert a(s)\vert
ds\right]^p+\frac{p}{1-p}\int_{t}^{T}\left[\frac{s}{\int_{s}^{+\infty}\tau^2\vert a(\tau)\vert d\tau}\right]^{1-p}s^2\vert a(s)\vert ds
\end{eqnarray*}
for all $T\geq t\geq t_0$.

The estimates
\begin{eqnarray*}
&&\frac{T^{1-p}}{1-p}\left[\int_{T}^{+\infty}s^2\vert a(s)\vert
ds\right]^p=\frac{T^{1-p}}{1-p}\int_{T}^{+\infty}\left[\frac{1}{\int_{T}^{+\infty}\tau^2\vert a(\tau)\vert d\tau}\right]^{1-p}s^2\vert a(s)\vert ds\\
&&\leq\frac{1}{1-p}\int_{T}^{+\infty}\left[\frac{s}{\int_{s}^{+\infty}\tau^2\vert a(\tau)\vert d\tau}\right]^{1-p}s^2\vert a(s)\vert ds
\end{eqnarray*}
allow us to establish that
\begin{eqnarray*}
\frac{1}{[k(c+\eta)]^p}\int_{t}^{T}\vert y_{0}(s)\vert^p ds\leq\frac{1+p}{1-p}\int_{t}^{+\infty}\left[\frac{s}{\int_{s}^{+\infty}\tau^2\vert
a(\tau)\vert d\tau}\right]^{1-p}s^2\vert a(s)\vert ds.
\end{eqnarray*}

The conclusion follows by letting $T\rightarrow+\infty$.

The proof is complete. $\square$

\begin{example}
\emph{An example of functional coefficient $a$ in the linear case that verifies the hypotheses of Theorem \ref{oscil_L_p} is given by the
formula
\begin{eqnarray*}
t^2a(t)=b(t)=\left\{
\begin{array}{l}
a_k(t-9k),\thinspace t\in[9k,9k+1],\\
a_k(9k+2-t),\thinspace t\in[9k+1,9k+3],\\
a_k(t-9k-4),\thinspace t\in[9k+3,9k+4],\\
a_k(9k+4-t),\thinspace t\in[9k+4,9k+5],\\
a_k(t-9k-6),\thinspace t\in[9k+5,9k+7],\\
a_k(9k+8-t),\thinspace t\in[9k+7,9k+8],\\
0,\thinspace t\in[9k+8,9(k+1)],
\end{array}
\right.\qquad k\geq1.
\end{eqnarray*}
Here, we take $a_k=k^{-\alpha}-(k+1)^{-\alpha}$ for a certain integer $\alpha>\frac{2-p}{p}$.}

\emph{To help the computations, the $k$-th "cell" of the function $b$ can be visualized next.}

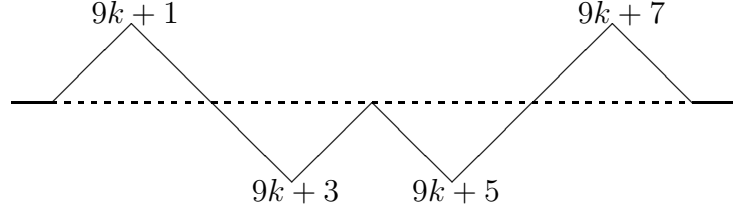
\begin{picture}(50,100)
\put(45,40){\line(1,0){15}} \put(60,40){\line(1,1){30}} \put(120,40){\line(-1,1){30}} \put(150,10){\line(-1,1){30}} \put(150,10){\line(1,1){30}}
\put(210,10){\line(-1,1){30}} \put(210,10){\line(1,1){60}} \put(300,40){\line(-1,1){30}} \put(300,40){\line(1,0){15}}
\multiput(60,40)(5,0){49}{\line(1,0){2}} \put(75,70){$9k+1$} \put(135,3){$9k+3$} \put(195,3){$9k+5$} \put(257,70){$9k+7$}
\end{picture}

\emph{It is easy to observe that
\begin{eqnarray*}
\int_{9k}^{9k+4}b(t)dt=\int_{9k+4}^{9k+8}b(t)dt=0\qquad\mbox{for all }k\geq1.
\end{eqnarray*}
}

\emph{We have
\begin{eqnarray*}
\int_{9k+2}^{+\infty}b(t)dt=\int_{9k+2}^{9k+4}b(t)dt=-a_k,\quad \int_{9k+6}^{+\infty}b(t)dt=\int_{9k+6}^{9k+8}b(t)dt=a_k
\end{eqnarray*}
and respectively
\begin{eqnarray*}
\int_{9k+2}^{+\infty}\vert b(t)\vert dt=3a_k+4\sum\limits_{m=k+1}^{+\infty}a_m,\qquad \int_{9k+6}^{+\infty}\vert b(t)\vert
dt=a_k+4\sum\limits_{m=k+1}^{+\infty}a_m.
\end{eqnarray*}
}

\emph{By noticing that
\begin{eqnarray*}
L_{+}=\lim\limits_{k\rightarrow+\infty}\frac{(9k+6)\int_{9k+6}^{+\infty}b(t)dt}{\int_{9k+6}^{+\infty}\vert b(t)\vert dt},\quad
L_{-}=\lim\limits_{k\rightarrow+\infty}\frac{(9k+2)\int_{9k+2}^{+\infty}b(t)dt}{\int_{9k+2}^{+\infty}\vert b(t)\vert dt},
\end{eqnarray*}
we obtain $L_{+}=\frac{9\alpha}{4}$ and $L_{-}=-\frac{9\alpha}{4}$.}

\emph{To verify the condition (\ref{hyp_a_L_p}), notice first that
\begin{eqnarray*}
I_{k}&=&\int_{9k}^{9(k+1)}\left[\frac{t}{\int_{t}^{+\infty}\vert b(s)\vert ds}\right]^{1-p}t^2\vert a(t)\vert
dt\\
&\leq&\int_{9k}^{9(k+1)}\left[\frac{9(k+1)}{\int_{9(k+1)}^{+\infty}\vert b(s)\vert ds}\right]^{1-p}a_{k}dt,\qquad k\geq1.
\end{eqnarray*}
}

\emph{The elementary inequality $a_k\leq(2^{\alpha}-1)(k+1)^{-\alpha}$ implies that
\begin{eqnarray*}
I_{k}\leq\frac{c_{\alpha}}{(k+1)^{(1+\alpha)p-1}},\qquad\mbox{where }c_{\alpha}=9\left(\frac{9}{4}\right)^{1-p}(2^{\alpha}-1),
\end{eqnarray*}
and the conclusion follows from the convergence of the series $\sum\limits_{k\geq1}(k+1)^{1-(1+\alpha)p}$.}
\end{example}

\textbf{Acknowledgement} The author is indebted to Professor Ondrej Do\v{s}ly and to a referee for valuable comments leading to an improvement
of the initial version of the manuscript. The author was financed during this research by the Romanian AT Grant 97GR/25.05.2007 with the CNCSIS
code 100.

\end{document}